\documentclass[11pt,oneside,a4paper]{amsart}
\usepackage{amssymb,latexsym}
\usepackage[british]{babel}
\usepackage{comment}
\usepackage{titletoc}
\usepackage{amsfonts, amsmath, amsthm, graphicx, color, times}
\usepackage[utf8]{inputenc}
\usepackage{enumerate}
\usepackage{xcolor}

\newtheorem {thm}{Theorem}
\newtheorem {cor}[thm]{Corollary}
\newtheorem {lem}[thm]{Lemma}
\newtheorem {prop}[thm]{Proposition}

\theoremstyle{definition}
\newtheorem {defi}[thm]{Definition}
\newtheorem {rem}[thm]{Remark}
\newtheorem {exa}[thm]{Example}

\DeclareMathOperator{\Gal}{Gal}

\DeclareMathOperator{\Ind}{Ind}

\DeclareMathOperator{\lcm}{lcm}
\DeclareMathOperator{\dens}{dens}
\DeclareMathOperator{\rk}{rank}

\newcommand{\p}{\mathfrak{p}}

\newcommand{\rank}{\text{rank}}

\newcommand{\Haar}{\mathrm{Haar}}
\newcommand{\mc}{\mathcal}

\usepackage{xr-hyper}
\usepackage{hyperref}

\title{Unified treatment of Artin-type problems II}

\begin{document}

\author{Olli~Järviniemi,  Antonella~Perucca and Pietro~Sgobba}
\address[]{Department of Mathematics and Statistics, University of Turku, 20014 Turku, Finland}
\email{olli.a.jarviniemi@utu.fi}
\address[]{Department of Mathematics, University of Luxembourg, 6 av.\@ de la Fonte, 4364 Esch-sur-Alzette, Luxembourg}
\email{antonella.perucca@uni.lu, pietrosgobba1@gmail.com}


\maketitle

\begin{abstract} 
This work concerns Artin's Conjecture on primitive roots and related problems for number fields. Let $K$ be a number field and let $W_1$ to $W_n$ be finitely generated subgroups of $K^\times$ of positive rank. We consider the index map, which maps a prime $\mathfrak p$ of $K$ to the $n$-tuple of the indices of  $(W_i \bmod \mathfrak p)$. Conditionally under GRH, any preimage under the index map admits a density, and the aim of this work is describing it. For example, we express the density as a limit in various ways. We study in particular the preimages of sets of $n$-tuples that are defined by prescribing valuations for their entries. Under some mild assumptions we can express the density as a multiple of a (suitably defined) Artin-type constant.
\end{abstract}


\section{Introduction}

Let $K$ be a number field, and work inside a fixed algebraic closure of $K$. Let $\alpha\in K^\times$ be not a root of unity, and write $\alpha^{1/m}$ for an $m$-th root of $\alpha$. Consider the primes $\mathfrak p$ of $K$ such that the reduction $(\alpha \bmod \p)$ is well-defined and non-zero, so that it makes sense to consider its index $\Ind_\p(\alpha)$ inside the multiplicative group of the residue field at $\mathfrak p$. 

Results on Artin's Conjecture for primitive roots by Hooley \cite{hooley}, Cooke and Weinberger \cite{cooke} (that are conditional under GRH)
ensure that the set of primes $\p$ such that $\Ind_\p(\alpha)=1$ has a density, which can be expressed in terms of the degrees of cyclotomic-Kummer extensions:
\begin{equation}\label{Hooley}
\dens\big\{\p:\Ind_\p(\alpha)=1\big\}=\sum_{n\geqslant 1}\frac{\mu(n)}{[K(\zeta_{n},\alpha^{1/n}):K]}\,.
\end{equation}
More generally, we consider a function $f(n)$ on the positive integers and the formal expression
\begin{equation}\label{formal}
\sum_{n\geqslant 1}\frac{\mu(n)}{[K(\zeta_{f(n)},\alpha^{1/f(n)}):K]}\,.
\end{equation}
Fixing some positive integer $t$, Ziegler \cite{ziegler}
proved (conditionally under GRH) that the set of primes $\p$ such that 
$\Ind_\p(\alpha)=t$ has a density given by \eqref{formal} setting $f(n)=nt$.
Moreover, Lenstra \cite{lenstra}
proved (conditionally under GRH) that the set of primes $\p$ such that 
$\Ind_\p(\alpha)\mid t$ has a density given by \eqref{formal} setting $f(n)=n\cdot\prod_{\ell\mid n} \ell^{v_\ell(t)}$.

For $k\geqslant 1$ we say that an integer is \emph{$k$-free} if it is not divisible by a $k$-th power greater than $1$: square-free is the same as $2$-free, and the number $1$ is the only $1$-free positive integer.
As a special case of the results in this paper, conditionally under GRH we have
\begin{equation}\label{kfree}
\dens\big\{\p:\Ind_\p(\alpha)\,\text{is $k$-free}\big\}=\sum_{n\geqslant 1}\frac{\mu(n)}{[K(\zeta_{n^k},\alpha^{1/{n^k}}):K]}\,.
\end{equation}

More generally, we could fix a set $\mathcal P'$ of prime numbers and for every $\ell\in \mathcal P'$ some positive integer $k(\ell)$: we may then require that the $\ell$-adic valuation of $\Ind_\p(\alpha)$ is less than  ${k(\ell)}$. Writing $\mathcal P'_\infty$ for the set of positive integers whose prime divisors all lie in $\mathcal P'$, conditionally under GRH we have 
\begin{equation}\label{formal2}
\dens\big\{\p: \;\forall \ell\in \mathcal P' \;\; \ell^{k(\ell)}\nmid \Ind_\p(\alpha)\big\} = \sum_{n\in \mathcal P'_\infty}\frac{\mu(n)}{[K(\zeta_{f(n)},\alpha^{1/f(n)}):K]}\,,
\end{equation}
where $f(n)=\prod_{\ell\mid n} \ell^{k(\ell)}$.
Notice that for some conditions on the index (for example, $k$-free for $k\geqslant 2$) one can get unconditional results because the involved cyclotomic-Kummer extensions contain a `large' cyclotomic extension (see the method by Pappalardi \cite{pappa}). Notice that the results stated above for $\alpha$ can be straight-forwardly generalized to a finitely generated subgroup of $K^\times$.

This paper should be seen as the follow up of \cite{JP} by the first two authors. Our results are again conditional under GRH, and we work in the same generality, namely we consider  groups $W_1,\ldots, W_n$ that are finitely generated subgroups of $K^\times$ of positive rank. We consider any set $H\subset \mathbb Z^n_{>0}$ and the density
\begin{equation}\label{density}
\dens\big\{\p: (\Ind_\p(W_1), \ldots, \Ind_\p(W_n))\in H\big\}\,,
\end{equation}
which is known to exist by \cite{JP}. 

We say that $H$ is \emph{cut by  valuations} if we have $H=\cap_\ell H_\ell$, where $H_\ell$ is the preimage under $v_\ell$ of $v_\ell(H)$: this  means that $H$ contains all integers that have suitable $\ell$-adic valuations for every $\ell$, the conditions on the various $\ell$ being independent. Slightly more generally, we allow finitely many valuation conditions not to be independent (we say that the set is \emph{almost cut by  valuations}). More generally, we consider a set $H$ which is \emph{determined by valuations}, by which we mean that $H=\cap_Q H_Q$, where $Q>1$ is square-free and $H_Q$ is the preimage of the $Q$-adic valuation of $H$.

Concerning the groups $W_1,\ldots, W_n$, we occasionally require that they are \emph{separated}, which means that no $W_i$ is contained in the group generated by the $W_j^{1/\infty}$ for $j\neq i$. This condition plays a role for the Kummer extensions and consequently for the index map 
$$\mathfrak p \mapsto (\Ind_\p(W_1), \ldots, \Ind_\p(W_n))\,,$$
see \cite{JP}. Notice that auxiliary results of Kummer theory of independent interest are proven in Section \ref{secKummer}.

General results about expressing \eqref{density} as a limit (in various ways) are contained in Section \ref{seclimits}: Proposition \ref{genlim} does not require any additional assumption, Theorem \ref{thmsep} only requires the groups to be separated, and Theorem \ref{thmdet} holds for all sets of tuples that are determined by valuations. 

Section \ref{secArtin} is devoted to define Artin-type constants that represent heuristical densities. By Proposition \ref{IMO} such constants are  strictly positive if the groups are separated. Then for separated groups and sets which are almost cut by valuations, we can express the density in \eqref{density} as a multiple of the Artin-type constant, see Remark \ref{multiple} and Proposition \ref{multiple2}. Finally, we can also handle correction factors, see Remark \ref{correction}.

\section{Cyclotomic-Kummer theory}\label{secKummer}

Let $K$ be a number field, and let $W_1$ to $W_n$ be finitely generated subgroups of $K^\times$ of positive rank. We write $I=\{1,\ldots, n\}$, and  
for every $J\subset I$ we write $W_J$ for the smallest subgroup of $K^\times$ containing the groups $W_i$ for all $i\in J$ (we set $W_\emptyset=\{1\}$).

We make use of the letter $\ell$ only to denote a prime number. 
Consider an $n$-tuple $e_I \in \mathbb{Z}_{\geqslant 0}^n$, denoting by $e_i$ its entries and by $e$ its maximum.
Calling $W_I^{1/\ell^{e_I}}$ the list $W_1^{1/\ell^{e_1}}, \ldots , W_n^{1/\ell^{e_n}}$, our aim is describing the degree
\begin{equation}
\label{eq:deg_kummer}
[K(\zeta_{\ell^\infty}, W_I^{1/\ell^{e_I}}) : K(\zeta_{\ell^\infty})]\,.    
\end{equation}
Up to reordering the groups, we will suppose that the tuple $e_I$ is \emph{non-increasing}, namely that $e_1 \geqslant \ldots \geqslant e_n$. Let $C$ be a constant and consider the partition $\mathcal{J}_C$ of $I$ into intervals, beginning a new interval at $i+1$ whenever $e_i - e_{i+1}> C$.
Let $J=[m,M]$ be an interval in $\mathcal{J}_C$. If $m\neq M$ (respectively, $m=M$), then we call \emph{difference tuple} of $e_I$ w.r.t.\ $C$ the $(M-m)$-tuple
\begin{equation}\label{differencetuple}
(e_m-e_{m+1}, \ldots, e_{M-1}-e_M)\in [0,C]^{M-m}\cap \mathbb{Z}_{\geqslant 0}^n
\end{equation}
(respectively, $\emptyset$).
Notice that by varying $e_I$ there are only finitely many possibilities for the difference tuple.

We call $\delta_J\in \mathbb Z^n_{\geqslant 0}$ the $n$-tuple such that $\delta_i=1$ for $i\in J$ and $\delta_i=0$ for $i\notin J$. 

\begin{lem}
\label{lem:inc_one}
For every sufficiently large $C$ (larger than a constant depending only on $K$ and $W_1,\ldots, W_n$), if $J=[m,M]\in \mathcal{J}_C$ is such that $m=1$ or $e_{m-1}-e_m>C+1$, then we have
$$v_\ell \big( [K(\zeta_{\ell^{\infty}}, W_I^{1/\ell^{e_I+\delta_J}}):K(\zeta_{\ell^{\infty}}, W_I^{1/\ell^{e_I}})] \big) = \begin{cases}
\rk(W_{[1,M]}) & \text{if $m=1$}\\
\rk(W_{[1,M]}) - \rk(W_{[1,m-1]}) & \text{if $m\neq 1$}\,.
\end{cases}$$
\end{lem}
\begin{proof}
We may assume w.l.o.g.\ that $W_I$ is torsion-free. Let $\mathcal{J_C}=\{J_1, \ldots, J_k\}$ and  write $J_i=[m_i, M_i]$. We define
$$\mathcal{W}_i = \prod_{j = m_i}^{M_i} W_{j}^{\ell^{e_{m_i} - e_j}}$$
thus we have
\begin{align*}
K(\zeta_{\ell^{\infty}}, W_I^{1/\ell^{e_I}}) = K(\zeta_{\ell^{\infty}}, \mathcal{W}_1^{1/\ell^{e_{m_1}}}, \ldots,  \mathcal{W}_k^{1/\ell^{e_{m_k}}}).
\end{align*}
Let $\rho_1=\rank(\mathcal{W}_1)$ and for each $i>1$ let $\rho_i = \rank(\mathcal{W}_{[1,i]}) - \rank(\mathcal{W}_{[1,i-1]})$. 
We inductively construct (for $i=1,\ldots, k$) multiplicatively independent elements $s_{i, j}$, where $1 \leqslant i \leqslant k$ and $1 \leqslant j \leqslant \rho_i$ , letting for each $i$ the elements $s_{i, 1}, \ldots , s_{i, \rho_i} \in \mathcal{W}_i$ be such that
$$\mathcal{W}_{i} \subset \langle \{ s_{i', j} | i' \leqslant i, j \leqslant \rho_{i'} \} \rangle^{1/\infty}.$$
We can choose the elements $s_{i, j}$ in a way that they depend only on the groups $W_1,\ldots, W_n$. Thus there is some positive integer $C_0$ (that depends only on the groups $W_1,\ldots, W_n$ because there are only finitely many possible values of $k$ and $(\mathcal{W}_1, \ldots , \mathcal{W}_k)$) such that 
$$\mathcal{W}_{i} \subset \langle \{ s_{i', j} | i' \leqslant i, j \leqslant \rho_{i'} \} \rangle^{1/C_0} \qquad \forall \,i=1,\ldots, k\,.$$

For every $C \geqslant v_{\ell}(C_0)$ we then have
\begin{align*}
K(\zeta_{\ell^{\infty}}, \mathcal{W}_1^{1/\ell^{e_{m_1}}}, \ldots , \mathcal{W}_k^{1/\ell^{e_{m_k}}}) = K(\zeta_{\ell^{\infty}}, \{s_{i, j}^{1/\ell^{e_{m_i}}} | 1 \leqslant i \leqslant k, 1 \leqslant j \leqslant \rho_i\}).
\end{align*}
By the same procedure we have
\begin{align*}
K(\zeta_{\ell^{\infty}}, W_I^{1/\ell^{e_I+\delta_J}}) = K(\zeta_{\ell^{\infty}}, \{s_{i, j}^{1/\ell^{e_{m_i}+\delta_{m_i}}} | 1 \leqslant i \leqslant k, 1 \leqslant j \leqslant \rho_i\}).
\end{align*}
so we conclude by Lemma \ref{Kummerlem}.
\end{proof}

\begin{lem}[{\cite[Proposition 3.1]{JP}}]\label{Kummerlem}
If $W$ is a finitely generated subgroup of $K^{\times}$, then there exists a positive integer $z$ (depending on $K$ and $W$, computable with an explicit finite procedure) such that for all positive integers $t$ we have
$$[K(\zeta_\infty, W^{1/zt}) : K(\zeta_\infty, W^{1/z})]=t^{\rk W}\,.$$
\end{lem}

We may now compute the degree \eqref{eq:deg_kummer} for any $\ell$ large enough:

\begin{lem}
\label{lem:large_ell}
There is a constant $c$ (depending only on $K$ and $W_1, \ldots , W_n$) such that for every $\ell > c$ and for every non-increasing $e_I\in \mathbb Z^n_{\geqslant 0}$ we have
\begin{align*}
v_{\ell}([K(\zeta_{\ell^{\infty}}, W_I^{1/\ell^{e_I}}) : K(\zeta_{\ell^{\infty}})]) = \sum_{1 \leqslant i \leqslant n} e_i(\rk(W_{[1,i]}) - \rk(W_{[1,i-1]}))\,.
\end{align*}
\end{lem}
\begin{proof}
Consider the proof of Lemma \ref{lem:inc_one}. If we fix some large enough $\ell$, then $v_{\ell}(C_0)=0$ and we may take $C=0$. We conclude because 
$$v_\ell \big([K(\zeta_{\ell^{\infty}}, \{s_{i,j}^{1/\ell^{e_{m_i}}} | 1 \leqslant i \leqslant k, 1 \leqslant j \leqslant \rho_i\}): K(\zeta_{\ell^{\infty}})]\big)=\sum_{i = 1}^k \rho_ie_{m_i}$$
and this sum equals the sum in the statement.
\end{proof}

Fix some $C$ such that Lemma \ref{lem:inc_one} applies. We call $\mathcal{T}(e_I)$ the difference tuple of $e_I$ w.r.t.\ $C$ and we call $\mathcal{T}$ the (finite) set of all difference tuples w.r.t.\ $C$ by varying $e_I$.
The next result describes the degree \eqref{eq:deg_kummer} for all $\ell$:
\begin{thm}
\label{thm:poly_degree}
There exist polynomials $P_{\ell, T}(x_I)\in \mathbb Z[x_I]$, $T \in \mathcal{T}$ such that the following holds: for every non-increasing $e_I\in \mathbb Z_{\geqslant 0}^n$ with $\mathcal{T}(e_I) = T$ we have
$${P_{\ell, T}(e_I)}=v_\ell\big([K(\zeta_{\ell^{\infty}}, W_I^{1/\ell^{e_I}}):K(\zeta_{\ell^{\infty}})]\big)\,.$$
There exist constants $c_{\ell, T}$ such that
\begin{align}
\label{eq:P}
P_{\ell, T}(x_I) + c_{\ell, T} = \sum_{i = 1}^n x_i\big(\rk(W_{[1,i]}) - \rk(W_{[1,i-1]})\big) \,.
\end{align}
Moreover, for every $\ell$ large enough we have $c_{\ell, T} = 0$ for every $T\in \mathcal T$.
\end{thm}
\begin{proof}
Fix $T$. Let $e_I'$ be the lexicographically smallest tuple such that $\mathcal{J}(e_I') = T$. Any other tuple $e_I$ such that $\mathcal{J}(e_I) = T$ is obtained by some sequence of operations as in Lemma \ref{lem:inc_one}, by first increasing the elements in the interval containing $e_1'$ a suitable number of times, then increasing the elements in the next interval and so on. Such an operation made for an interval $J=[m,M]$ increases the $\ell$-adic valuation of the degree of the relevant extension by $\rank(W_{[1,M]}) - \rank(W_{[1,m]})$. It follows that the requested degree is of form $\ell^{P_{\ell, T}(e_I)}$ for some polynomial $P_{\ell, T}$ as in \eqref{eq:P}. The fact that $c_{\ell, T} = 0$ for large $\ell$ is Lemma \ref{lem:large_ell}. 
\end{proof}

\begin{rem}
In Theorem \ref{thm:poly_degree} we have $c_{\ell, T}\geqslant 0$ for any multiplicatively independent groups because $\rank(W_i)=\rank(W_{[1,i]}) - \rank(W_{[1,i-1]})$ and hence, letting $e'_I$ be as in the proof of the theorem, we have
$$v_\ell \big([K(\zeta_{\ell^{\infty}}, W_I^{1/\ell^{e'_I}}):K(\zeta_{\ell^{\infty}})]\big)\leqslant \sum_{1 \leqslant i \leqslant n} e'_i(\rank(W_{[1,i]}) - \rank(W_{[1,i-1]})\big)\,.    
$$
\end{rem}

\section{Limits of densities}\label{seclimits}

\subsection{Valuation conditions}
We consider the \emph{$\ell$-adic valuation} as a map
$v_\ell: \mathbb Z^n_{>0}\rightarrow \mathbb  Z^n_{\geqslant 0}$ and more generally for a square-free integer $Q>1$ we define the \emph{$Q$-adic valuation} $v_Q$ as the product of the maps $v_\ell$ for $\ell\mid Q$. Let $H\subset \mathbb Z^n_{>0}$, and write $H_Q$ for the preimage under $v_Q$ of $v_Q (H)$. 
For every $\ell$ large enough we have $0_I\in H_\ell$.
For all square-free integers $Q,Q_0>1$ with $Q_0\mid Q$ we have 
\begin{equation}\label{sets}
H\subset H_Q \subset H_{Q_0} \cap \bigcap_{\ell\mid \frac{Q}{Q_0}} H_\ell\subset \bigcap_{\ell \mid Q} H_{\ell}\,.
\end{equation}

\begin{defi}\label{defisets}
We say that $H$ is \emph{cut by valuations} if $H=\cap_{\ell} H_\ell$.
We say that $H$ is \emph{almost cut by valuations} if
there is some square-free integer $Q_0>1$ such that
$$H=H_{Q_0} \cap \bigcap_{\ell\nmid Q_0} H_\ell\,.$$
We say that $H$ is \emph{determined by valuations} if $H=\cap_Q H_Q$ by varying $Q>1$ square-free.
\end{defi}

By \eqref{sets} a set which is cut by valuations is almost cut by valuations, and a set which is almost cut by valuations is determined by valuations.

\begin{exa}
The set of positive square-free integers is cut by valuations. The set of positive integers with an even amount of prime divisors is not determined by valuations because we have  $H_Q=\mathbb Z$ for every $Q$.
A finite set or, more generally, a set with finite support (the set of prime divisors of the entries is finite), is not necessarily cut by valuations, but it is almost cut by valuations.  
An example of set that is determined by  valuations but it is not almost cut by valuations is the set of prime numbers.
\end{exa}

\subsection{Notation}

We keep the notation of \cite{JP}. In particular, $\ell$ is a prime number and $Q \geqslant 1$ is square-free. For $h_I\in \mathbb Z_{>0}^n$ we write $h:=\lcm(h_1,\ldots, h_n)$.
Moreover, we write $h_{Q,I}$ (respectively, $h_Q$) for the truncation of $h_I$ (respectively, $h$) obtained by removing the prime factors coprime to $Q$.

We fix some finite Galois extension $F/K$ and a union of conjugacy classes $C\subset  \Gal(F/K)$.
We write 
 $$K_{Q^\infty,Q^\infty_I}=\bigcup_{\ell\mid Q, e\geqslant 1} K(\zeta_{\ell^e}, W_I^{1/\ell^e})\,.$$
For $H\subset \mathbb Z_{>0}^n$ we define 
$$C_{H,Q}:=\bigcup_{h_I\in H} C_{h_{Q, I}}=\mathop{\dot{\bigcup}}_ {h_{Q,I}} C_{h_{Q, I}}\subset \Gal(FK_{Q^\infty,Q^\infty_I}/K)\,,$$
where $C_{h_{Q, I}}$ is the conjugacy-stable set of those automorphisms $\sigma$ that satisfy all of the following conditions: the restriction of $\sigma$ to $F$ lies in $C$; the restriction of $\sigma$ to $K(\zeta_{h_{Q}}, W_I^{1/h_{Q,I}})$
is the identity; 
for all $q\mid Q$ and for all $i\in I$ the restriction of $\sigma$ to $K(\zeta_{qh_{Q,i}}, W_i^{1/qh_{Q,i}})$ is not the identity.

Consider the Haar measure $\mu_{\Haar}$ on  $\Gal(FK_{Q^\infty,Q^\infty_I}/K)$. The set $C_{h_{Q, I}}$ (respectively, $C_{H,Q}$) is measurable and it has measure $0$ if and only if it is empty. Moreover, we have 
$$\mu_{\Haar}(C_{H,Q})=\sum_{h_{Q,I}} \mu_{\Haar}(C_{h_{Q, I}})\,.$$

\subsection{Results}
Consider the primes $\mathfrak p$ of $K$ unramified in $F$ for which the index map 
\begin{equation}\label{indexmap}
\Psi: \mathfrak p \rightarrow (\Ind_{\p}W_1,\ldots, \Ind_{\p}W_n)
\end{equation}
is well-defined.
Let the set $S_{C,H}$ consist of those $\mathfrak p$ such that  $\binom{\p}{F/K}\subset C$ and $\Psi(\p)\in H$.

\begin{thm}[{\cite{JP}}]\label{thmJP}
Assume GRH. The set $S_{C,H}$ admits a natural density, and we have
\begin{equation}\label{eq-singleton}
\dens (S_{C,H})=\sum_{h_I\in H} \dens (S_{C,\{h_I\}})\,.
\end{equation}
Moreover, we have $\dens (S_{C,H})=0$ if and only if $S_{C,H}=\emptyset$.
\end{thm}

The sum in the above statement means considering the $n$-tuples  such that the maximum of the entries is at most $x$, and letting $x\rightarrow \infty$. In the statements below the limit for $Q\rightarrow \infty$ means the limit for $Q_x:=\prod_{\ell\leq x}\ell$ while the limits $Q,B\rightarrow \infty$ and $Q,k\rightarrow \infty$ can be seen as double limits (where the order of the limits does not matter) but also as a single limit, for example $\min(Q,B)\rightarrow \infty$.

\begin{prop}\label{genlim}
Assume GRH. Then we have 
\begin{equation}\label{limB}
\dens(S_{C,H})=\lim_{B\rightarrow \infty} \dens (S_{C,H\cap [1,B]^n})\,.
\end{equation}
If $Q$ is a positive square-free number, calling $Q^\infty$ the subset of $\mathbb Z_{>0}^n$ consisting of the $n$-tuples such that the prime divisors of each of the entries divide $Q$, we have 
\begin{equation}\label{limQinfty}
    \dens(S_{C,H})= \lim_{Q\rightarrow \infty} \dens(S_{C,H\cap Q^\infty})
\end{equation}
and
\begin{equation}\label{limQinftyB}
    \dens(S_{C,H})=
 \lim_{Q,B\rightarrow \infty} \dens(S_{C,H\cap Q^\infty\cap [1,B]^n})
    \,.
\end{equation}
Calling $\mathbb Z^n_ {\text{$k$-free}}$ the tuples consisting of $k$-free integers, we also have
\begin{equation}\label{limQinftykfree}
\dens(S_{C,H}) = \lim_{Q,k\rightarrow \infty} \dens(S_{C,H\cap Q^\infty\cap \mathbb Z^n_ {\text{$k$-free}}})\,.
\end{equation}
\end{prop}
\begin{proof}
By Theorem \ref{thmJP} all densities exist. Since $[1,B]^n\subset [1,B']^n$ for $B\leq B'$, the limit in \eqref{limB} exists and the inequality $\geqslant$ holds. Thus to prove \eqref{limB} it  suffices to show  that the the density of the set $\{\p : \Psi(\p)\notin [1,B]^n\}$ goes to $0$ for $B\rightarrow \infty$.
This set is contained in the finite union for $i=1,\ldots,n$ of the set $\{\p : \Ind_\p (W_i)>B \}$
and we conclude because 
$$\sum_{t\geqslant 1} \dens(\{\mathfrak p : \Ind_\p (W_i)=t\})=1<\infty\,.$$
Equality \eqref{limQinfty} follows from $H=\bigcup_Q (H\cap Q^\infty)$ and $\sum_{h\in H} \dens(S_{C,\{h\}})<\infty$.
For \eqref{limQinftyB} consider that  $[1,B]^n\subset Q^\infty$ holds for some suitable $Q$ (namely, the square-free part of $B!$). Finally, we can deduce \eqref{limQinftykfree} from \eqref{limQinftyB} because $[1,B]^n\subset \mathbb Z^n_ {\text{$k$-free}}$ holds for some suitable $k$ (namely, the largest valuation of the numbers in $[1,B]$). \end{proof}

\begin{thm}\label{thmsep}
Assume GRH, and suppose that the groups are separated. Then the following are equivalent:
\begin{itemize}
\item[(i)] $\dens (S_{C,H}) >0$;
\item[(ii)] $C_{H,Q} \neq \emptyset$ holds for every square-free integer $Q>1$;
\item[(iii)] $C_{H,Q_x} \neq \emptyset$ holds for some sufficiently large integer $x$ (which  is computable and depends on $K$ and $W_1,\ldots, W_n$ but it does not depend on $H$).
\end{itemize}
\end{thm}
\begin{proof}
The implications (i)$\Rightarrow$(ii)$\Rightarrow$(iii) are clear.
The implication (iii)$\Rightarrow$(i) holds for a singleton by \cite[Theorem 4.1]{JP}. In general, as the groups are separated, we can make use of 
\cite[Theorem 6.1 (iii)]{JP}.
\end{proof}

The conditions in the above theorem are clearly necessary to have $\dens (S_{C,H}) >0$. They are also sufficient if $H$ is a singleton by \cite{JP} but in general they are not: the following example can be  generalized to any groups that are not separated, and it shows that we cannot have (ii)$\Rightarrow$(i).

\begin{exa}\label{cexa}
Consider $n=2$, $K=F=\mathbb Q$, $W_1 = W_2 = \langle 2 \rangle$ and 
$H = \{(q, q^2) |\, \text{$q$ prime}\}$.
Clearly $S_{C,H}=\emptyset$.
However, $v_Q(0_I)\in v_Q(H)$ holds for every $Q$ hence $C_{H,Q}\neq \emptyset$ (by considering only $v_Q$ we are not excluding those primes $p$ such that $2\bmod p$ is a primitive root).
\end{exa}

\begin{thm}\label{thmdet}
Assume GRH and suppose that $H$ is determined by valuations. Then we have 
\begin{equation}\label{Haar2}
 \dens (S) = \lim_{Q \to \infty} \mu_{\Haar} (C_{H,Q})
 \end{equation}
 and 
\begin{equation}\label{limit2}
 \dens (S) = \lim_{Q \to \infty}\dens\bigg\{\p: \Psi(\p)\in H_Q,\ \begin{pmatrix}\p \\ F/K  \end{pmatrix} \subset C\bigg\}\,.
 \end{equation}
\end{thm}
\begin{proof}
Formula \eqref{Haar2} is a consequence of \eqref{limit2}.
To prove \eqref{limit2} it suffices to show that
the density of the set $\big\{\p: \Psi(\p)\in H_Q\setminus H \big\}$ goes to $0$ for $Q\rightarrow \infty$.
Call $H'=\mathbb Z_{>0}^n \setminus H$. To every $t_I\in H'$ we can associate the least positive integer $z$ such that $t_I\notin H_{Q_x}$ whenever $x\geqslant z$.
This function provides a partition of $H'$ into sets $H'_z$. Moreover, by Theorem \ref{thmJP} we have 
$$\dens(H')=\sum_{z\geqslant 1} \dens(H'_z)\,.$$
As this series converges, we conclude because the complement of $H_{Q_x}$ is contained in  $\cup_{z\geqslant x} H'_z$.
\end{proof}

In the following statement we let $Q_0$ be as in Definition \ref{defisets} (notice that in the limit we could restrict to those $\ell\mid Q$ such that $\ell\nmid Q_0$).

\begin{cor}\label{corcuts}
Assume GRH and suppose that $H$ is almost cut by valuations. Then we have 
 \begin{equation}\label{limit1}
 \dens (S) = \lim_{Q \to \infty}\dens\bigg\{\p: \Psi(\p)\in H_{Q_0} \cap \bigcap_{\ell\mid Q}H_\ell,\ \begin{pmatrix}\p \\ F/K  \end{pmatrix} \subset C\bigg\}\,.
 \end{equation}
If $H$ is cut by valuations, then we have 
 \begin{equation}\label{limit1}
 \dens (S) = \lim_{Q \to \infty}\dens\bigg\{\p: \Psi(\p)\in \bigcap_{\ell\mid Q}H_\ell,\ \begin{pmatrix}\p \\ F/K  \end{pmatrix} \subset C\bigg\}\,.
 \end{equation} 
\end{cor}
\begin{proof}
This is a consequence of Theorem \ref{thmdet} by writing $H_Q$ suitably.
\end{proof}

\begin{exa}
We cannot expect Theorem \ref{thmdet} to hold for all sets $H$. For example, let $n=1$, consider a trivial Frobenius condition, and let $H$ be the set of prime numbers. Then $\cap_\ell H_\ell$ is the set of all positive square-free integers and we have $1\in \cap_Q H_Q$.
So by considering $\cap_\ell H_\ell$ or $\cap_Q H_Q$ we may obtain a strictly larger density (provided that  there is a positive density of primes $\mathfrak p$ such that $\Ind_\p(W_1)=1$).
\end{exa}

\section{Artin-type constants}\label{secArtin}

To study the density of the preimage of a set $H\subset \mathbb Z^n_{>0}$ under the index map \eqref{indexmap}, we define an  Artin-type constant, whose value depends on $H$ and the ranks of the groups $W_J, J \subset I$.

For every $i$, let $r_i:=\rank(W_i)$. Also let $\mathcal{R}$ be the tuple of integers $R_J := \rank(W_J), J \subset I$. We assume $R_J > 0$ holds for $J \neq \emptyset$, and we set $R_\emptyset=0$. In case for every $J$ we have $R_J = \sum_{j \in J} r_j$ (corresponding to multiplicatively independent groups),  we write $r_I$ instead of $\mathcal{R}$.

We define the \emph{inclusion-exclusion function} on a  tuple of real numbers $x_I$ as follows (for the empty tuple we have $\mc {IE}(x_I)=1$):
\begin{equation}
\mathcal I\mathcal E(x_I):=\sum_{ J\subset I}(-1)^{\#J} \prod_{j\in J}x_j=\prod_{i\in I} (1-x_i)\,. \end{equation}
Notice that setting $x_i=\frac{1}{d_i}$ the function $\mathcal I\mathcal E$ provides the density of the primes of a number field that do not split in any of finitely many linearly disjoint Galois extension of the number field having degree $d_i$.

For any $\ell$ we let $V_\ell\subset \mathbb Z^n_{\geqslant 0}$, and we write $V:=\prod_\ell V_\ell$. Given $H$, we may set $V_\ell:=v_\ell(H)$ and moreover, every $V$ such that $0_I\in V_\ell$ holds for every $\ell\gg 1$ can be obtained in this way.
We define an Artin-type constant 
\begin{equation}
\label{eq:artin_constant}
\mathcal A_{V,\mathcal R}:=\prod_{\ell} A_{V_\ell,\mathcal R}(\ell)
\end{equation}
as an infinite product (which converges) of suitably defined real numbers in  $[0,1]$.
In turn, we define  
\begin{equation}
A_{V_\ell,\mathcal R}(\ell):=\sum_{v_I\in V_\ell} F_{v_I,\mathcal R}(\ell)
\end{equation}
as a series (which converges) of suitably defined real numbers in $[0,1]$. The numbers $F_{v_I,\mathcal  R}(\ell)$ are the evaluation at $\ell$ of a rational function $F_{v_I,\mathcal R}(x)$ with integer coefficients, corresponding to a heuristical density of primes $\mathfrak{p}$ so that $v_{\ell}(\Ind_{\mathfrak{p}}(W_i)) = v_i \text{ for all } i \in I.$

Let us define these heuristical densities $F_{v_I, \mathcal{R}}$. We first note that the event $v_{\ell}(\Ind_{\mathfrak{p}}(W)) \geqslant 1$ ought to happen with probability $\ell^{-\rank(W)}$ assuming that $\mathfrak{p}$ splits in $K(\zeta_{\ell})$. This leads to the following formula in the simple case of the zero $n$-tuple $0_I$.

\begin{defi}
For the zero $n$-tuple, we define
\begin{equation}\label{Artin-zerotuple}
    F_{0_I,\mc R}(\ell) :=
\frac{\ell-2}{\ell-1}+\frac{1}{\ell-1}\sum_{J\subset I}\frac{(-1)^{\# J}}{\ell^{R_J}} = 1 - \frac{1}{\ell-1} \sum_{\emptyset \neq J \subset I} \frac{(-1)^{\# J}}{\ell^{R_J}}\,.
\end{equation}
\end{defi}

Now consider an $n$-tuple $v_I\neq0_I$ and call $v:=\max_i(v_i)$. There is a permutation $\sigma$ such that the tuple $(v_{\sigma_1}, \ldots, v_{\sigma_n})$ is non-increasing. Then we define a function $f:\mathbb Z_{\geqslant 0}^n\rightarrow\mathbb Z$ by setting
$$f(x_I)= x_{\sigma_1}R_{\sigma_1} + \sum_{i=2}^n x_{\sigma_i} \big(R_{\{\sigma_1,\ldots,\sigma_i\}}-R_{\{\sigma_1,\ldots,\sigma_{i-1}\}}\big)\,.$$
Notice that the choice of $\sigma$ does not affect the value of $f$ (to see this it suffices to make swaps of neighboring indices ${\sigma_i}, {\sigma_{i+1}}$ such that $x_{\sigma_i}=x_{\sigma_{i+1}}$).

\begin{defi}
For $v_I\neq0_I$ we define
\begin{equation}  \label{eq-Fgeneral}
\begin{split}
F_{v_I, \mathcal R}(\ell) & := \frac{1}{\varphi(\ell^v)\ell^{f(v_I)}} \left(\frac{\ell - 1}{\ell} \sum_{\substack{J \subset I \\ J \cap I' = \emptyset}} \frac{(-1)^{|J|}}{\ell^{f(v_I + \delta_J) - f(v_I)}} + \frac{1}{\ell} \sum_{\substack{J \subset I}} \frac{(-1)^{|J|}}{\ell^{f(v_I + \delta_J) - f(v_I)}}\right)\,.
\end{split}
\end{equation}
\end{defi}
Let us motivate the definition. First, the factor $1/\varphi(\ell^v)$ comes from the density of $\mathfrak{p}$ splitting in $K(\zeta_{\ell^v})$. The factor $1/\ell^{f(v_I)}$ comes from requiring that $\mathfrak{p}$ splits in $K(\zeta_{\ell^v}, W_i^{1/\ell^v})$ for all $i$: conditionally on $\mathfrak{p}$ splitting in 
$$K(\zeta_{\ell^v}, W_{\sigma_1}^{1/\ell^{v_{\sigma_1}}}, \ldots , W_{\sigma_{i-1}}^{1/\ell^{v_{\sigma_{i-1}}}}),$$
there is a probability 
$$\ell^{-v_{\sigma_i}(R_{\{\sigma_1, \ldots , \sigma_i\}} - R_{\{\sigma_1, \ldots , \sigma_{i-1}\}})}$$
that $\mathfrak{p}$ splits in
$$K(\zeta_{\ell^v}, W_{\sigma_1}^{1/\ell^{v_{\sigma_1}}}, \ldots , W_{\sigma_{i}}^{1/\ell^{v_{\sigma_{i}}}})$$
as well. We then consider the probability that $\mathfrak{p}$ does not split in any of $K(\zeta_{\ell^{v_i + 1}}, W_i^{1/\ell^{v_i + 1}})$, considering separately two cases depending on whether $\mathfrak{p}$ splits in $K(\zeta_{\ell^{v+1}})$. The first sum in \eqref{eq-Fgeneral} corresponds to the case where it does not, while the second sum corresponds to the case where it does, the formulas being obtained by inclusion-exclusion and similar reasoning as above.

While it is not obvious from the definition \eqref{eq-Fgeneral}, the values $F_{v_I, \mathcal{R}}(\ell)$ are non-negative and under a separateness condition strictly positive, see Lemma \ref{IMO}. The formula may be rewritten as
\begin{equation}
F_{v_I, \mathcal{R}}(\ell) = \frac{1}{\ell^v}\Bigg(\sum_{J\subset I \atop J\cap I'=\varnothing} \frac{(-1)^{\# J}}{\ell^{f(v_I+\delta_J)}} + \frac{1}{\ell-1} \sum_{J\subset I}  \frac{(-1)^{\# J}}{\ell^{f(v_I+\delta_J)}}\Bigg)\,.
\end{equation}

\begin{exa}\label{exa1}
For $n=1$ (setting $r=r_1$) we have
$$F_{v,r}(\ell)=\begin{cases}
1-\frac{1}{\ell^r(\ell-1)}& \text{for $v=0$}\\
\frac{1}{\ell^{v(r+1)}}\cdot \frac{\ell}{(\ell-1)} (1-\frac{1}{\ell^{r+1}}) & \text{for $v>0$\,.}
\end{cases}
$$
\end{exa}

\begin{exa}
For a constant $n$-tuple $v_I\neq 0_I$ we have
\[ F_{v_I,\mc R}(\ell) = \frac{1}{\ell^{v(1+R_I)}}\cdot\Big(1+\frac{1}{\ell-1}\sum_{J\subset I}\frac{(-1)^{\# J}}{\ell^{R_J}}\Big) \]
because we can write $f(v_I+\delta_J) = vR_I+R_J$.
\end{exa}

\begin{exa}\label{exa2}
For $n=2$, supposing that $R_I=r_1+r_2$, we have
\[ F_{v_I,r_I}(\ell)=\begin{cases}
1+\frac{1}{\ell-1}(\mc {IE}(\ell^{-r_1}, \ell^{-r_2})-1)& \text{for $v_1=v_2=0$}\\
\frac{1}{\ell^{v(1+r_1+r_2)}}\cdot \frac{\ell}{\ell-1} \cdot \big(1+\frac{1}{\ell}\cdot (\mc {IE}(\ell^{-r_1}, \ell^{-r_2})-1)\big) & \text{for $v_1= v_2>0$}\\
\frac{1}{\ell^{v+v_1r_1+v_2r_2}}\cdot  
\frac{\ell}{\ell-1}\cdot 
\mathcal I\mathcal E(\ell^{-r_m}, \ell^{-(r_M+1)})
& \text{for $v_1\neq v_2$}\\
\end{cases}
\]
where $I=\{m,M\}$ is such that $v_M=v$.
\end{exa}

\begin{exa}\label{ex-independent}
Suppose that for every $J$ we have $R_J=\sum_{j\in J} r_j$. Then we have
$$F_{0_I,r_I}(\ell):=\frac{\ell-2}{\ell-1}+ \frac{1}{\ell-1}\mc{IE}(\ell^{-r_I}) \,.$$
For $v_I\neq0$ we have
\begin{align}
F_{v_I,r_I}(\ell) 
&  = \frac{1}{\ell^{v+\sum_iv_ir_i}}  \cdot\left( \mc{IE}(\ell^{-r_{I\setminus I'}})+\frac{1}{\ell-1}\mc{IE}(\ell^{-r_{I}})  \right)  \label{eq-Fv}\\
& = \frac{1}{\ell^{v+\sum_iv_ir_i}} \cdot \frac{\ell}{\ell-1} \cdot \mc{IE}(\ell^{-r_{I\setminus I'}})
\cdot \Big( 1+\frac{1}{\ell} \cdot  ( \mc{IE}(\ell^{-r_{I'}})  -1) \Big)\notag
\end{align}
because we can write $f(v_I+\delta_J) - \sum_{i\in I} v_ir_i =\sum_{i\in J} r_i$.
\end{exa}

We then note that the probabilistic motivation for \eqref{eq-Fgeneral} to be the density of the preimage is rigorous for $\ell$ sufficiently large.

\begin{cor}\label{cor-ell-large}
There exists a constant $c$ such that for every $\ell>c$ the preimage of $v_I$ under $v_\ell\circ \Psi$ has density $F_{v_I, \mathcal R}(\ell)$.
\end{cor}
\begin{proof}
Suppose that $\ell\gg 1$.
For $v_I=0_I$, we have
$$    F_{0_I,\mc R}(\ell) =
\frac{\ell-2}{\ell-1}+\frac{1}{\ell-1}\sum_{J\subset I}\frac{(-1)^{\# J}}{\ell^{R_J}}\,.$$
The first summand is the density of primes not splitting completely in $K(\zeta_{\ell})$, while $1/\ell$ is the density of those that do. Then we are left to remark that $[K(\zeta_\ell, W_J^{1/\ell}):K(\zeta_\ell)]=\ell^{-R_J}$.
For $v_I\neq 0_I$ we have 
$$F_{v_I, \mathcal R}(\ell) =   
\frac{1}{\ell^v}\Bigg(\sum_{J\subset I \atop J\cap I'=\varnothing} \frac{(-1)^{\# J}}{\ell^{f(v_I+\delta_J)}} + \frac{1}{\ell-1} \sum_{J\subset I}  \frac{(-1)^{\# J}}{\ell^{f(v_I+\delta_J)}}\Bigg).$$
The first summand accounts for the suitable primes splitting completely in $K(\zeta_{\ell^v})$ but not in $K(\zeta_{\ell^{v+1}})$, while the second summand accounts for the suitable primes splitting completely in $K(\zeta_{\ell^{v+1}})$. We are left to prove that for $J\subset I\setminus I'$ we have (by varying $i\in I\setminus J$ and $j\in J$)
$$[K(\zeta_{\ell^v}, W_{i}^{1/\ell^{v_i}}, W_{j}^{1/\ell^{v_j+1}}):K(\zeta_{\ell^v})]=\ell^{-f(v_I+\delta_J)}$$
and that for all $J$ we have 
$$[K(\zeta_{\ell^{v+1}}, W_{i}^{1/\ell^{v_i}}, W_{j}^{1/\ell^{v_j+1}}):K(\zeta_{\ell^{v+1}})]=\ell^{-f(v_I+\delta_J)}\,.$$
Up to reordering the groups, we may suppose that the tuple $w_I:=v_I+\delta_J$ is non-increasing and hence 
$$f(w_I)=
\sum_{i\in I}w_{i} \big(R_{[1,i]}-R_{[1,i-1]}\big)\,.$$
We conclude because for every $M\geqslant \max_i w_i$ and for every $i\in I\setminus\{1\}$ adding the radicals $W_i^{1/\ell^{w_i}}$ to the field 
$$K(\zeta_{\ell^M}, W_1^{1/\ell^{w_1}}, \ldots, W_{i-1}^{1/\ell^{w_{i-1}}})$$ is the same as adding the radicals ${\mathcal W}_i^{1/\ell^{w_i}}$, where ${\mathcal W}_i$ is a torsion-free subgroup of $W_i$ of rank $R_{[1,i]}-R_{[1,i-1]}$ such that the groups $W_i$ and ${\mathcal W}_i$ have the same rank modulo $W_{[1,i-1]}$ (this is because the tuple $w_I$ is non-increasing and because the radicals of the torsion elements in that quotient produce a trivial extension).
\end{proof}

Finally, we note that if the groups $W_1, \ldots , W_n$ are separated, then $F_{v_I, \mathcal{R}}(\ell)$ indeed is strictly positive, as is the resulting Artin-type constant $A_{V, \mathcal{R}}$.

\begin{prop}\label{IMO}
Suppose that $W_1, \ldots , W_n$ are separated. Then for every $\ell$ and for every $v_I\in \mathbb Z^n_{\geqslant 0}$ we have $F_{v_I, \mathcal R}(\ell)>0$. Moreover, if $0_I\in V_{\ell}$ holds for $\ell\gg 1$, then the Artin-type constant $A_{V,\mathcal R}$ is strictly positive.
\end{prop}
\begin{proof}
The second claim follows from the first because for $\ell\gg 1$, supposing w.l.o.g.\ that $V_\ell=\{0_I\}$ holds, we have 
$$F_{0_I,\mathcal R}(\ell)\geqslant 1 -\frac{1}{\ell-1}\sum_{\varnothing\neq J\subset I} \frac{1}{\ell}>1-\frac{2^n}{\ell^2-\ell}\,.$$

The claim on the positivity of $F_{v_I, \mathcal{R}}(\ell)$ is proven via a probabilistic interpretation. We may assume w.l.o.g.\ that $W_I$ is torsion-free. For $1 \leqslant i \leqslant n$, let $b_{i, 1}, \ldots,  b_{i, r_i}$ be a basis for $W_i$ and let $B$ denote the multiset $\{b_{i, j} | 1 \leqslant i \leqslant n, 1 \leqslant j \leqslant r_i\}$. Associate to each $b_{i, j}$ a random variable $X_{i, j}$, for which
$$P(X_{i, j} = k) = \frac{1}{\ell^k} - \frac{1}{\ell^{k+1}} \text{\; for } 0 \leqslant k < v+1, \quad P(X_{i, j} = v+1) = \frac{1}{\ell^{v+1}},$$
with the random variables $X_{i, j}$ independent of each other. Introduce a further random variable $X_{\zeta}$ with
$$P(X_{\zeta} = k) = \frac{1}{\varphi(\ell^k)} - \frac{1}{\varphi(\ell^{k+1})} \text{\; for } 0 \leqslant k < v+1, \quad P(X_{\zeta} = v+1) = \frac{1}{\varphi(\ell^{v+1})},$$
with $X_{i, j}$ and $X_{\zeta}$ all independent of each other. (Intuitively $X_{\zeta}$ and $X_{i, j}$ correspond to the largest $v$ such that a random prime $\mathfrak{p}$ splits in $K(\zeta_{\ell^k})$ or $K(\zeta_{\ell^k}, b_{i, j}^{1/\ell^k})$.)

Let $E$ denote the event
$$\text{for any } B' \subset B \text{ and } b_{i, j} \in B \setminus B' \text{ we have } b_{i, j} \in \langle B' \rangle^{1/\infty} \implies X_{i, j} \geqslant \min_{b_{i', j'} \in B'} X_{i', j'}.$$
Note that $P(E) \geqslant P(X_{i, j} = v+1 \text{ for all } i, j) > 0$ hence, we may condition on $E$. (The event $E$ accounts for the fact that the groups $W_i$ and fields $K(\zeta_{\ell^k}, W_i^{1/\ell^k})$ may interact with each other.)

Consider then the conditional probability
$$P(\text{for all } i \in [1, n], \min (X_{i, 1}, \ldots , X_{i, r_i}, X_{\zeta}) = v_i | E).$$
The probabilistic motivation given for \eqref{eq-Fgeneral} shows that this probability equals $F_{v_I, \mathcal{R}}(\ell)$. This already proves that $F_{v_I, \mathcal{R}}(\ell)$ is non-negative. 
For strict positivity, note that by the assumption on the separateness of $W_1, \ldots , W_n$, one may construct indices $t_1, \ldots , t_n$ such that for any $i$, the random variable $X_{i, t_i}$ is independent from 
$$\{X_{i', j'} | 1 \leqslant i' \leqslant n, j' \leqslant r_{i'}\} \setminus \{X_{i, t_i}\}$$
conditionally on $E$. Hence
\begin{align*}
&P(\text{for all } i \in [1, n], \min (X_{i, 1}, \ldots , X_{i, r_i}, X_{\zeta}) = v_i | E) \geqslant \\
&P(X_{\zeta} = v+1 \text{ and for all } i \in [1, n], X_{i, t_i} = v_i \text{ and for all } j \neq t_i, X_{i, j} = v+1 | E) > 0.
\end{align*}
\end{proof}

\section{Euler products and multiplicative correction factors}\label{secEuler}

In the following statement $Q_0$ is as in Definition \ref{defisets}.

\begin{prop}\label{almostcut}
There is some square-free integer $B\geqslant 1$ (depending only on $K$ and $W_1,\ldots,W_n$) such that if $H$ is almost cut by valuations, then we have
\begin{equation}\label{eq-densH}
\dens(S_{C,H}) = \mu_{\Haar}(C_{H_{\lcm(B,Q_0)}}) \cdot \prod_{\ell\nmid \lcm(B,Q_0)}\mu_{\Haar}(C_{H_\ell}).
\end{equation}
\end{prop}
\begin{proof}
By Lemma \ref{Kummerlem} there is a square-free integer $B$ such that, calling $B':=\lcm(B,Q_0)$ and considering square-free integers $Q>1$ such that $B'\mid Q$, we have  
\[ C_{H,Q}=C_{H_B'}\times\bigoplus_{\ell\nmid B',\ell\mid Q}C_{H_\ell}. \]
We then conclude by \eqref{Haar2}.
\end{proof}

\begin{rem}\label{multiple}
For every $\ell\gg 1$ we have $\mu_{\Haar}(C_{H_\ell})=A_{V_\ell,\mc R}$. So if $W_1,\ldots, W_n$ are separated and $H$ is almost cut by valuations, then by Propositions \ref{IMO} and \ref{almostcut} we obtain 
\begin{equation}\label{Euler}
 \dens(S_{C,H}) = \mc A_{V,\mc R}\cdot \bigg( \mu_{\Haar}(C_{H_B})\cdot\prod_{\ell\mid B} A_{V_\ell,\mc R}^{-1}\bigg)\,.
\end{equation}
We deduce that $\dens(S_{C,H})$ is strictly positive if $\mu_{\Haar}(C_{H_B})>0$.
Notice that, in case $C_{H_B}$ is the preimage of its projection to a finite Galois group, $\dens(S_{C,H})$ is by \eqref{Euler} a rational multiple of Artin's constant $\mc A_{V,\mc R}$.
\end{rem}

\begin{rem}\label{correction}
For general sets $H$, one method for computing $\dens (S_{C,H})$ is provided by \eqref{eq-singleton}. Indeed, by \cite[Theorem 4.1]{JP} for every $h_I\in \mathbb Z_{>0}^n$ we know that $\dens (S_{C,\{h_I\}})$ is a rational multiple of the strictly positive Artin-type constant $\mathcal A_{V, \mathcal R}$ with $V_\ell=\{0_I\}$ for every $\ell$ (which amounts to the heuristical density for the index to be $1$). Calling $m(h_I)$ the multiplicative correction factor, we get
\begin{equation}\label{formulagen}
\dens (S_{C,H})=A_{V, \mathcal R} \cdot \sum_{h_I\in H} m(h_I)\,.
\end{equation}
For some ``easy'' sets $H$ the sum over correction factors may be simplified via other means.
\end{rem}

Notice that, for separated groups, altering $h_I$ with prime factors $\ell\gg 1$ does not affect the positivity of the density. In fact, we have the following:

\begin{prop}\label{multiple2}
Suppose that $W_1,\ldots, W_n$ are separated. There is some square-free integer $Q\geqslant 1$ such that for every $h_I\in \mathbb Z^n_{>0}$ we have 
$$\dens (S_{C,\{h_I\}}) = \mu_{\Haar}(C_{\{h_I\},Q}) \cdot \prod_{\ell\nmid Q} F_{v_\ell(h_I), \mathcal R}(\ell)\,.$$
In particular, $\dens (S_{C,\{h_I\}})$ is strictly positive if $\mu_{\Haar}(C_{\{h_I\},Q})>0$.
\end{prop}
\begin{proof}
We may apply \eqref{Haar2}. Since the groups are separated there is some square-free integer $B\geqslant 1$ (independent of $h_I$) such that for any square-free integer $Q\geqslant 1$ we have
$$\mu_{\Haar}(C_{\{h_I\}, Q}) = \mu_{\Haar}(C_{\{h_I\}, \gcd(B, Q)}) \prod_{\ell \mid Q, \ell \nmid B} f_{\ell, \{h_I\}},$$
where, calling $h_{\ell,i}:=\ell^{v_\ell(h_i)}$, and $h_\ell:=\max_I h_{\ell,i}$,  the rational number $f_{\ell, \{h_I\}}\in [0,1]$ is the proportion of automorphisms of 
$$K(\zeta_{\ell h_\ell}, W_1^{1/\ell h_{\ell, 1}}, \ldots , W_n^{1/\ell h_{\ell, n}})/K(\zeta_{h_\ell}, W_1^{1/h_{\ell, 1}}, \ldots , W_n^{1/h_{\ell, n}})$$
that do not restrict to the identity on any of the fields $K(\zeta_{\ell h_{\ell, i}}, W_i^{1/\ell h_{\ell, i}})$. Assuming $B$ is large enough, this proportion is by definition equal to $F_{v_\ell(h_I), \mathcal{R}}(\ell)$.
\end{proof}

\begin{exa}
Let $n=1$, $\rank(W_1)=1$, $K=F=\mathbb Q$ and let $H$ be the set of prime numbers. For example, consider the density of the primes $p$ such that the index of $(2 \bmod p)$ is prime.
Call $A_{t}$ the density of the primes $p$ for which $\Ind_p W_1=t$. Suppose that $A_{1}>0$, which means that the generator of $W_1$ is not a perfect square. Then for all but finitely many primes $q$ we have 
by Example \ref{exa1}
$$\frac{A_q}{A_1}=
\bigg(\frac{1}{q(q-1)} \Big(1-\frac{1}{q^{2}}\Big)\bigg)\cdot 
\bigg(1-\frac{1}{q(q-1)}\bigg)^{-1}=\frac{q^2-1}{q^2 (q^2-q-1)}\,.$$
In the generic case, the requested density is then 
$$A \cdot \sum_q \frac{q^2-1}{q^2 (q^2-q-1)}\,,$$
where $A$ is the usual Artin's constant. 
In the non-generic case the density is of the form 
$$c A + c' A_1 \sum_q \frac{q^2-1}{q^2 (q^2-q-1)}$$
where $c,c'$ are rational numbers and here $A_1$ is also a rational multiple of $A$.
Notice that if $A_1=0$ (which means that the generator of $W_1$ is a perfect square) the index $\Ind_p(W_1)$ cannot be an odd prime, so it is prime if and only if it equals $2$.
\end{exa}

\end{document}